\documentstyle[12pt,fleqn,leqno]{article}
\newtheorem{thm}{Theorem}[section]
\newcommand{\be}{\begin{equation}}
\newcommand{\ee}{\end{equation}}
\newcommand{\ba}{\begin{array}}
\newcommand{\ea}{\end{array}}

\renewcommand{\a}{\alpha}
\renewcommand{\b}{\beta}
\renewcommand{\l}{\lambda}

\renewcommand{\em}{\it}
\newcommand{\bea}{\begin{eqnarray}}
\newcommand{\eea}{\end{eqnarray}}

\begin{document}
\newtheorem{lem}[thm]{Lemma}
\newtheorem{cor}[thm]{Corollary}
\title{\bf Ladder Operators For Szeg\H{o} Polynomials and Related Biorthogonal Rational Functions\thanks{
Research partially supported  by NSF grant DMS 9203659 and NSERC grant A6197} }
\author{ Mourad E. H.  Ismail and Mizan Rahman}
\date{}

\maketitle
\begin{abstract} We find  the raising and lowering operators for
 orthogonal polynomials on the unit circle introduced by Szeg\H{o}
 and for their four parameter generalization to ${}_4\phi_3$ biorthogonal
 rational functions
 on the unit circle.
\end{abstract}

\bigskip 
{\bf Running title}: Ladder Operators.



\bigskip

\setcounter{section}{1}

\setcounter{equation}{0}

{\bf 1. Introduction}. The Szeg\H{o} theory for orthogonal polynomials on
 the unit circle \cite{Gr:Sz} was
 developed in the early part of  this century. One aspect of this theory is
 to relate the limiting
 behavior of the orthonormal polynomials to the absolutely continuous of part of
 the measure with respect to which the polynomials are orthonormal. 
As an illustration of his theory Szeg\H{o} \cite{Sz} proved that the polynomials
 $\{{\cal H}_n(z|q)\}$,
\be
{\cal H}_n(z|q) := \sum_{k=0}^n \frac{(q;q)_n}{(q;q)_k(q;q)_{n-k}}
 (q^{-1/2}z)^k,
\ee
satisfy the orthogonality relation
\be
\frac{1}{2\pi i}\int_{|z| = 1}\overline{{\cal H}_m(z|q)}{\cal H}_n(z|q)
w_c(z|q) \frac{dz}{z}
 = q^{-n}\frac{(q;q)_n}
{(q;q)_\infty} \delta_{m,n},
\ee
where the weight function $w_c(z|q)$ is given by
\be
w_c(z|q) := (q^{1/2}z, q^{1/2}/z; q)_\infty.
\ee
Throughout this paper we shall always assume that $0 < q < 1$ and shall adopt the following 
notation for the $q$-shifted factorials \cite{Ga:Ra}
\be
(a; q)_0 := 1, \quad (a; q)_n := \prod_{k = 1}^n (1 - a q^{k-1}),\quad n =
1, 2, \cdots, \; or\;  \infty,
\ee
and for the multiple $q$-shifted factorials 
\be
(a_1, a_2, \cdots, a_k)_n := \prod_{j = 1}^k (a_j; q)_n.
\ee
We shall also need the notion of basic hypergeometric series
\bea
{}_{r}\phi_s\left(\left. \ba{c} a_1,\ldots,a_{r} \\ b_1,\ldots,b_s \ea \right|
\,q, \;z \right) &=& {}_{r}\phi_s(a_1, \ldots, a_{r}; b_1, 
\ldots, b_s; q, z) \\
&=& \sum_{n=0}^\infty 
\frac{(a_1,\ldots,a_{r};q)_n}{(q,b_1,\ldots,b_s;q)_n}\,z^n (-q^{n(n-1)/2})^{s
+1-r}. \nonumber
\eea
The subscript $c$ in $w_c(z|q)$ refers to the unit circle. 

There is a hierarchy of  the so called classical $q$-orthogonal polynomials
 with the Askey-Wilson polynomials (\cite {An:As}, \cite{As:Wi}, \cite{Ga:Ra})
 being the most general and the continuous
 $q$-Hermite polynomials (\cite{As:Is1}, \cite{Ga:Ra}) at the bottom of the chart.
The polynomials at a certain level are special or limiting cases of  some of the 
polynomials at
 a higher level. In \cite{Be:Is} Berg and Ismail showed how one can use generating
 functions to go from the continuous $q$-Hermite polynomials to the Askey-Wilson
  polynomials  through the intermediate stage of the Al-Salam-Chihara
 polynomials, \cite{Al:Ch}.  We believe that the Szeg\H{o} polynomials $\{{\cal H}_n(z|q)\}$ are
 the  unit circle version
of the continuous  $q$-Hermite polynomials. The unit circle analogues of the
 Askey-Wilson polynomials are the four-parameter biorthogonal 
 rational functions introduced by Al-Salam and Ismail in \cite{Al:Is}. They are the
 pair $\{r_n\}$, $\{s_n\}$,
\bea
r_n(z; a, \a, b, \b|q) := {}_4\phi_3\left(\left. \ba{c} q^{-n},
 ab\a\b q^{n-1}, bq^{1/2},  b z\\
  b\a, b\b, abq^{1/2}z \ea \right|\,q,
\;q \right)  
\eea
and 
\be
 s_n(z;, a, \a, b, \b|q)  = r_n(z; \overline{\a},\overline{a},\overline{\b},
\overline{b}|q).
\ee
Al-Salam and Ismail \cite{Al:Is} established the biorthogonality relation
\bea
 &&\frac{1}{2\pi i} \int_{|z|=1} w_c(z;  a, \a, b, \b|q) 
 r_n(z; a, \a, b, \b|q)
\overline{s_m(z;  a, \a, b, \b|q)}\; \frac{dz}{z} \\
&\mbox{}& \quad = \kappa(a, \a, b, \b)\frac{(q, a\a, ab\a\b q^{n-1};q)_n}
{(b\b;q)_{n}(ab\a\b;q)_{2n}}
 (b\b)^n \delta_{m,n}, \nonumber
\eea
where the weight function $w_c(z;  a, \a, b, \b|q)  $ is
\bea
w_c(z;  a, \a, b, \b|q)  := \frac{(q^{1/2}z, q^{1/2}/z, abq^{1/2}z, \a\b q^{1/2}/z;q)_\infty}
{(az, \a /z, bz, \b/z;q)_\infty}.
\eea
and the total mass  $\kappa(a, \a, b, \b)$ is
\bea
\kappa(a, \a, b, \b) &= &  \frac{1}{2\pi i} \int_{|z| = 1} w_c(z;  a, \a, b, \b|q) 
\frac{dz}{z}\\
&=&
\frac{(aq^{1/2},\a q^{1/2}, bq^{1/2}, 
\b q^{1/2}, ab\a\b ;q)_\infty}{ (q, a\a, b\a, a\b, b\b;q)_\infty}. \nonumber
\eea
When $a = \a = 0$ the functions $r_n$ become polynomials previously studied 
by Pastro, \cite{Pa}.  
Observe that (1.11) is a special case of \cite[(4.11.3)]{Ga:Ra}.

The $q$-difference operator $D_{q,z}$ is 
\be
(D_{q,z}f)(z) := \frac{f(z) - f(qz)}{(1-q)z}.
\ee
The purpose of this work is to first compute the adjoint of $D_{q,z}$ on a 
suitable inner product space, then to show that $D_{q,z}$ and its adjoint are, 
respectively,  the lowering and  
raising operators of the Szeg\H{o} polynomials
 $\{{\cal H}_n(z|q)\}$ polynomials, and finally to generalize these results
 to the biorthogonal rational functions $\{r_n\}$ and $\{s_n\}$. In  Section 2 we shall treat the polynomials $\{{\cal H}_n(z|q)\}$ and we also find the
 $q$-Sturm-Liouville system whose eigenfunctions are the Szeg\H{o} polynomials.
 This latter fact is then
 used to give an elementary proof of the orthogonality of the ${\cal H}_n$'s
 with respect  to $w_c(z|q)$.
 Section 3 develops the $q$-Sturm-Liouville problem associated with $D_q$ on the unit circle. 
In Section  4 we identify the raising and lowering operators for the more general
 and more complicated case of the biorthogonal pair $\{r_n, s_n\}$. Section 5 contains
 a new derivation of the biorthogonality relation (1.9) through the use of raising and lowering
 operators found in Section 4.

This work is part of the renewed interest in biorthogonal rational functions 
and polynomials, \cite{Al:Ve}, \cite{Is:Ma}, \cite{Ra}, \cite{Ra:Su}. 

\bigskip

\setcounter{equation}{0}

\setcounter{section}{2}

\setcounter{thm}{0}

{2. \bf Ladder Operators}. The operator $D_{q,z}$ of (1.12) acts nicely on 
the ${\cal H}_n$'s. It is straightforward to derive
\be
D_{q,z} {\cal H}_n(z|q) = \frac{q^{-1/2}(1-q^n)}{1-q} {\cal H}_{n-1}(z|q),
\ee
from the representation (1.1).
This shows that $D_{q,z}$ acts as a lowering  operator on the ${\cal H}_n$'s. In order to
 find a  raising operator we need to compute the adjoint of $D_{q,z}$ with 
respect to a suitable inner product. Before we introduce the appropriate
 inner product space we need to 
introduce some notations. If $f$ is analytic in $\rho_1 < |z| < \rho_2$ then
$\overline{f}$ will denote the function whose Laurent coefficients  are the
 complex conjugates of
 the corresponding Laurent coefficients of $f$. Thus $\overline{f}$ is also 
analytic in
 $\rho_1 < |z| < \rho_2$. 
Consider the set of functions
\be
{\cal F}_\nu := \{f: f(z) \mbox{ is analytic for }   q^\nu \le |z| \le
q^{-\nu} \}.
\ee
This set will be equipped with the inner product  
\be
<f,g>_c := \frac{1}{2\pi i} \int_{|z| = 1} f(z) \overline{g(z)}\, \frac{dz}
{z}.
\ee
 We shall now use ${\cal F}_{\nu}$ to denote this inner product space. It is clear that
 if $f \in {\cal F}_{\nu}$ then $\overline{f} \in {\cal F}_{\nu}$.
\begin{thm}
The adjoint of the $q$-difference operator $D_{q,z}$  is 
$T_{q,z}$,
\be
(T_{q,z}f)(z)= \frac{z[f(z) - qf(qz)]}{1-q},
\ee
that is 
\be
<D_{q,z}f, g>_c = <f, T_{q,z}g>_c, \quad \mbox{for}\;  f, g \in {\cal F}_1.
\ee
\end{thm}
{\bf Proof}. We have
\bea
 <D_{q,z}f,g>_c &=& \frac{1}{2\pi i} \int_{|z|=1} \frac{f(z) - f(qz)}
{(1-q)z}\overline{g(z)}
\, \frac{dz}{z} \nonumber \\
&=&\frac{1}{2\pi i} \int_{|z|= 1} \frac{f(z)}{(1-q)z}
\overline{g}(z^{-1})\, \frac{dz}{z} - \frac{1}{2\pi i} \int_{|z|=q^{-1}}
 \frac{f(qz)}{(1-q)z}\overline{g}(z^{-1}) \frac{dz}{z}, \nonumber
\eea
since $f$ and $\overline{g}$ are analytic in $1\le |z| \le q^{-1}$. We 
replace $z$ by $q^{-1}z$ in the last integral. This gives
\bea
<D_{q,z}f,g>_c &=& \frac{1}{2\pi i} \int_{|z|=1} f(z)
 \frac{\overline{g}(z^{-1})-q\overline{g}(qz^{-1})}{(1-q)z}
\frac{dz}{z}\\
 &=& \frac{1}{2\pi i} \int_{|z|=1} f(z)
 \frac{\overline{z[g(z) -qg(qz)]}}{(1-q)}\frac{dz}{z}, \nonumber
\eea
and the proof is complete.

Observe that $T_{q,z}$ can be written in the form
\bea
(T_{q,z}f)(z) = qz^2 (D_{q,z}f)(z) + zf(z).
\eea

We now show that $T_{q,z}$ is a raising operator for the
 polynomials ${\cal H}_n(z|q)$. Since $w_c(z|q)$ is real on the unit circle
 we can rewrite the orthogonality relation (1.2)  in the form
\bea
q^{-n}\frac{(q;q)_n}{(q;q)_\infty}
 \delta_{m,n} &=& <{\cal H}_m(z|q),
 w_c(z|q){\cal H}_n(z|q)>_c \nonumber \\
&=& \frac{q^{1/2}(1-q)}{1-q^{m+1}}<D_{q,z} {\cal H}_{m+1}(z|q), 
w_c(z|q){\cal H}_n(z|q)>_c
 \nonumber \\
&=& \frac{q^{1/2}(1-q)}{1-q^{m+1}}< {\cal H}_{m+1}(z|q),
 T_{q,z} w_c(z|q){\cal H}_n(z|q)>_c.
 \nonumber
\eea
This shows that  
\bea
\frac{1}{w_c(z|q)}T_{q,z} w_c(z|q){\cal H}_n(z|q) \nonumber
\eea
is orthogonal to ${\cal H}_{n+1}(z|q)$ for all $m \ne n$.  Hence one would expect
\bea
\frac{1}{w_c(z|q)}T_{q,z} w_c(z|q){\cal H}_n(z|q) \nonumber
\eea
to be a constant multiple of ${\cal H}_{n+1}(z|q)$.
 Thus we have have been led to the following result.

\begin{thm} The raising operator for $\{{\cal H}_{n}\}$ is
 $T_{q,z}$ in the sense
\bea
\frac{1}{w_c(z|q)}T_{q,z} \left( w_c(z|q){\cal H}_n(z|q)\right)
 = \frac{\sqrt{q}}{1-q} \, {\cal H}_{n+1}(z|q).
\eea
\end{thm}

This theorem  follows by direct evaluation of the left-hand side of
 (2.8). The details are straightforward and are omitted. 

\begin{thm}
The Szeg\H{o} polynomials have the Rodrigues formula
\be
{\cal H}_n(z|q) =(q^{-1/2} - q^{1/2})^n \frac{1}{w_c(z|q)}T_{q,z}^n
 \left(w_c(z|q)\right).
\ee
\end{thm}
{\bf Proof}. Apply (2.8) repeatedly.

When we combine (2.4) and (2.8) we arrive at the following theorem.
\begin{thm}
The polynomials $\{{\cal H}_n(z|q)\}$ satisfy the $q$-Sturm-Liouville equation
\be
\frac{1}{w_c(z|q)}T_{q,z} \left(w_c(z|q)D_{q,z} {\cal H}_n(z|q)\right)
 = \l_n {\cal H}_n(z|q),
\ee
where 
\be
\l_n = \frac{(1-q^n)}{(1-q)^2}.
\ee
\end{thm}

Observe that the eigenvalues of (2.11) are distinct. 

\begin{thm} 
The   orthogonality relation (1.2) follows from (2.1) and (2.8).
\end{thm}
{\bf Proof}. Set
\bea
\zeta_n = <w_c(z|q) {\cal H}_n(z|q), {\cal H}_n(z|q)>. 
\eea
Clearly we have
\bea
\frac{q^{-1/2}(1-q^{n+1})}{1-q}\zeta_n &=& < D_{q,z}{\cal H}_{n+1}(z|q), 
w_c(z|q){\cal H}_n(z|q)>  \nonumber \\
&=& < {\cal H}_{n+1}(z|q), T_{q,z} w_c(z|q){\cal H}_n(z|q)> \nonumber\\
&=& \frac{q^{1/2}}{1-q} \zeta_{n+1}. \nonumber
\eea
Therefore $\zeta_{n+1} = q^{-1}(1-q^{n+1})\zeta_n$ and we get
\bea
\zeta_n = q^{-n}(q;q)_n\zeta_0.   
\eea
To find $\zeta_0$ expand $w_c(z|q)$ using the Jacobi triple product identity \cite[(II.28)]{Ga:Ra}
\bea
\sum_{-\infty}^\infty q^{n^2}z^n = (q^2, -qz, -q/z;q)_\infty.
\eea
It is then straightforward to see that $\zeta_0 = 1/(q;q)_\infty$. Now 
substitute for $\zeta_0$ in (2.13) to obtain (1.2) and this completes the proof.

\bigskip

\setcounter{equation}{0}

\setcounter{section}{3}

\setcounter{thm}{0}

{3. \bf $q$-Sturm Liouville Operators}. Motivated by equation (2.9)
we  now consider the more general operator
\be
(M f)(z) : = \frac{1}{\omega(z)}\left(T_{q,z} \left(p(z)D_{q,z} f\right) \right)(z),
\ee
where $\omega(z)$ is real on the unit circle with some restrictions on $p$ 
and $\omega$ to 
follow.  The analysis in this section follows closely the plan in \cite{Br:Ev:Is}, 
where the spectral theory of the Askey-Wilson polynomials is discussed. 
 
Let ${\cal H}_\omega$ denote the inner product space $L^2$ of the unit circle
equipped with the inner product
\be
(f, g)_\omega := \frac{1}{2\pi i}\int_{|z| = 1} f(z)\overline{g(z)} \omega(z) 
\frac{dz}{z},
\ee
and let
\bea
T := M_{|{\cal F}_2} \; in \; {\cal H}_w. \nonumber
\eea
We shall assume that $p$ and $\omega$ satisfy
\bea
&&(i)\quad \; p(z) > 0 \;a.e.\; on \; |z| = 1, \; p \in {\cal F}_1, \; 1/p \in
 L(\{z: |z| =1\}); \\
&&(ii)  \quad On \; the \; unit \; circle\; \omega(z) > 0\; a.e. \;and \; both
\;   \omega \; and \; 
 1/\omega \; are \; integrable. \nonumber
\eea
The expression $Mf$ is therefore defined for $f \in {\cal F}_2$, and the operator 
$T$ acts in ${\cal H}_w$. Furthermore, the domain ${\cal F}_2$ of $T$ is dense 
in ${\cal H}_\omega$ since it contains all polynomials and Laurent polynomials.
\begin{thm}
The operator $T$ is symmetric in ${\cal H}_\omega$ and $T \ge 0$.
\end{thm}
{\bf Proof}. For all $f, g \in {\cal F}_1$,
\bea
(f, Tf)_\omega &=&  <f, T_{q,z}\left( p(z)D_{q,z} f\right)>
 = <D_{q,z}f, p(z)D_{q,z} f)
\\
&=& \frac{1}{2\pi i}\int_{|z| = 1} p(z) |D_{q,z}f(z)|^2 \frac{dz}{z}, \nonumber
\eea
which proves our theorem.
\begin{cor}
Let $y_1,y_2 \in {\cal F}_2$ be solutions to 
\bea
\frac{1}{\omega(z)}\left(T_{q,z} \left(p(z)D_{q,z} f\right) \right)(z) = \l f,
\eea
 with $\l = \l_1$ and $\l = \l_2$, respectively,  and assume $\l_1\ne \l_2$. Then $y_1$
 and $y_2$ are orthogonal in the sense
\bea
\frac{1}{2\pi i}\int_{|z| = 1} y_1(z) y_2(z) \frac{dz}{z} = 0.
\eea
Furthermore, the eigenvalues of (3.5) are all real.
\end{cor}

Let ${\cal Q}(T)$ denote the form domain of $T$ and let $\tilde{T}$  be its 
Friedrichs extension. Recall that the domain of $T$, ${\cal Q}(T)$, is the 
completion of ${\cal F}_1$ with respect to $\|.\|_{{\cal Q}}$, where 
 \be
\|.\|^2_{{\cal Q}} := \frac{1}{2\pi i}\int_{|z| = 1} p(z) |D_{q,z} .|^2 
\frac{dz}{z} + \|f\|^2_\omega,
\ee
and if $(.,.)_{{\cal Q}}$ denotes the inner product on ${\cal Q}(T)$, then
 for all $f \in {\cal D}(\tilde{T})$ and $g \in {\cal Q}(T)$,
\be
(f, g)_{\cal Q} = (\tilde{T}f, g)_\omega.
\ee 
We have that $f \in {\cal Q}(T)$ if and only if there exists a sequence 
$\{f_n\}
\subset {\cal F}_1$ such that $\|f - f_n\|_{\cal Q} \to 0$; hence $\|f - f_n\|_w
 \to 0$ and $\{{\cal D}_qf_n\}$ is a Cauchy sequence in $L^2(\{z: |z| = 1\};
 p(z)\frac{dz}{2 \pi i z})$,
 with limit $F$ say. From (3.3) and (2.3) it follows that for
 $\phi \in {\cal F}_{1}$,
\bea
\frac{1}{2\pi i}\int_{|z| = 1} F(z) \phi(z) \frac{dz}{z} &=& 
\lim_{n \to \infty} \frac{1}{2\pi i}\int_{|z| = 1} (D_{q,z}f_n)
(z) \phi(z) \frac{dz}{z} \nonumber \\
 &=&  \lim_{n \to \infty} \frac{1}{2\pi i}\int_{|z| = 1} f_n(z) D_{q,z}
\phi(z) \frac{dz}{z},  \nonumber \\
   &=&  \frac{1}{2\pi i}\int_{|z| = 1} f(z) D_{q,z} \phi(z) \frac{dz}{z}.
\nonumber
\eea
Thus in analogy with distributional derivatives, we shall say that
 $F = D_{q,z} f$ in the generalized sense. We conclude that the 
norm on ${\cal Q}(T)$ is defined by (3.7) with $D_{q,z} f$ now 
defined in the generalized sense. Also, it follows in a standard way that
\bea
{\cal D}(T^*) = \{f: f, Mf \in {\cal H}_\omega\}, \quad Tf = Mf,
\eea
\bea
{\cal D}(\tilde{T}) =&& {\cal Q}(T) \cap {\cal D}(T^*)  \\
=&& \{f: p^{1/2} D_{q,z} f \in L^2(\{z: |z| =1\}), \; Mf \in {\cal H}_\omega\}. 
\nonumber
\eea

\bigskip

\setcounter{equation}{0}

\setcounter{section}{4}

\setcounter{thm}{0}

{4. \bf Biorthogonal Functions}. We first evaluate the action of
 $D_{q,z}$ on $r_n(z; a, \a, b, \b|q)$. A calculation gives
\bea
D_{q,z} \frac{(bz; q)_k}{(abq^{1/2}z;q)_k} 
= -\frac{b(1-aq^{1/2})(1-q^k)}{1-q} \;
\frac{(qbz; q)_{k-1}}{(abq^{1/2}z;q)_{k+1}}. \nonumber
\eea
Therefore
\bea
(abq^{1/2}z;q)_2 \,D_{q,z} r_n(z; a, \a, b, \b|q) &=& 
\frac{bq^{1-n}(1-aq^{1/2})(1-bq^{1/2})(1-q^{n})}
{(1-q)(1-b\a)(1- b\b)} \\
&\mbox{}& \quad \times (1-ab\a\b q^{n-1}) r_{n-1}(z; qa, \a, bq, \b|q). \nonumber
\eea
This shows that the operator $(abq^{1/2}z; q)_2 D_{q,z}$ is a lowering  operator for the
 $r_n$'s. To guess at the raising operator we use (2.5). 
 Clearly
\bea
&&<(abq^{1/2}z;q)_2 D_{q,z} f, w_c(1/z;q\overline{a}, \overline{\a},
 \overline{qb}, \overline{\b}|q) \; g> \nonumber \\
&& \mbox{} \quad \mbox{} \quad = <D_{q,z} f,(\overline{a}\overline{b}q^{1/2}/z;q)_2
 w_c(1/z; q\overline{a},
 \overline{\a}, \overline{qb}, \overline{\b}|q) \;g> \nonumber \\
&& \mbox{} \quad \mbox{} \quad =< f, T_{q,z} 
\left((\overline{a}\overline{b}q^{1/2}/z;q)_2
 w_c(1/z; q\overline{a},
 \overline{\a}, \overline{qb}, \overline{\b}|q))\; g\right)>. \nonumber
\eea
This relationship, (4.1)  and the orthogonality relation (1.9) show that when $m \ne n$
 then 
\bea
0 &=& <r_{m-1}(z; qa, \a, qb, \b|q) w_{c}(z; qa, \a, qb, \b|q),
s_{n-1}(z; qa, \a, qb, \b|q)>  \\
&=& <r_{n}(z; a, \a, b, \b|q), \nonumber \\
&& \mbox{} \quad T_{q,z}\left((\overline{a}\overline{b}q^{1/2}/z;q)_2
 w_c(1/z;q\overline{a}, \overline{\a},
 q\overline{b}, \overline{\b}|q) s_{n-1}(z; qa, \a, qb, \b|q)\right)>.
\nonumber
\eea
But (1.10) shows the symmetry of the weight function
\bea
w_c(1/z; a, \a, b, \b|q) = w_c(z; a, \a, b, \b|q).
\eea
Now (4.2), (4.3), (1.8) and the biorthogonality relation (1.9) suggest that
\bea
T_{q,z}\left((\overline{a}\overline{b}q^{1/2}/z;q)_2
 w_c(z;q\overline{a}, \overline{\a},
 q\overline{b}, \overline{\b}|q) s_{n-1}(z; qa, \a, qb, \b|q)\right)
\nonumber
\eea
is a multiple of
\bea
w_c(1/z; \overline{a} , \overline{\a}, \overline{b}, 
 \overline{\b}|q) s_{n}(z; qa, \a, qb, \b|q).
\eea
This suggests that the raising operator we are looking for is defined by
\be
(L^+f)(z) = \frac{1}{w_c(z; a, \a, b, \b|q)}T_{q,z}
 \left((\a\b q^{-3/2}/z;q)_2
 w_c(z; a, \a, b, \b|q)f(z)\right).
\ee
Observe that $L^+$ of (4.5) is a raising operator when it acts on $r_n$ 
because the functions $\{r_n, s_n\}$ are biorthogonal and not necessarily
 orthogonal. Like the case of the Szeg\H{o} polynomials 
 it turns out that we have made the correct guess.

\begin{thm} The raising and lowering operators for $\{r_n(z; a, \a, b, \b|q)\}$ 
are $L^+$ and $L^-$, where
\bea
(L^-f)(z): = ((\a\b q^{1/2}z; q)_2 D_{q,z}f)(z),
\eea
and $L^+$ is as in (4.5). In other words 
\bea
L^- r_{n}(z; a, \a, b, \b|q) &=& 
\frac{bq^{1-n}(1-q^{1/2}a)(1-q^{1/2}b)(1-ab\a\b q^{n-1})}
{(1-q)(1-b\a)(1-b\b)} \\
&\mbox{}& \quad \times (1-q^n) r_{n-1}(z; qa, \a,  qb, \b|q), \nonumber
\eea
and
\bea
L^+ r_{n-1}(z; a, \a, b, \b|q) = \frac{(1-b\a/q)(1-b\b/q)}{(1-q)b} 
 r_{n}(z;a, \a/q, b, \b/q|q)
\eea
\end{thm}
{\bf Proof}. We need only only to prove (4.8). Set
\bea
 && T_{q,z}  \left( (\a\b q^{1/2}/z;q)_2 w_c(z; a, q\a, b, q\b|q)
r_{n-1}(z; a, q\a, b, q\b|q) \right) \\
&& \mbox{} \quad = \frac{z}{1-q} w_c(z; a, \a, b, \b|q) \; L_n(z). \nonumber
\eea
A calculation using (4.9), (1.7) and (1.10) gives 
\bea
&&L_n(z) = (1- \a/z)(1-\b/z)
{}_4\phi_3\left(\left. \ba{c} q^{1-n},
 ab\a\b q^{n}, q^{1/2}b, bz\\
 qb\a, qb\b, abq^{1/2}z \ea \right|\,q,
\;q \right)  \\
&& \mbox{} \quad + \frac{(zq^{1/2} - \a\b )(1-az)(1-bz)}
{z^2\; (1- ab q^{1/2}z)}
{}_4\phi_3\left(\left. \ba{c} q^{1-n},
 ab\a\b q^{n}, q^{1/2}b,  qb z\\
 ab q^{3/2}z, qb\a, qb\b \ea \right|\,q,
\;q \right).  \nonumber
\eea
Recall the Sears transformation \cite[(III.15)]{Ga:Ra}
\bea
&& {}_4\phi_3\left(\left. \ba{c} q^{-n},
 A, B, C\\
  D, E, F \ea \right|\,q,
\;q \right)  \\
&& \mbox{} \quad = \frac{(E/A, F/A;q)_n}{(E, F;q)_n}A^n
{}_4\phi_3\left(\left. \ba{c} q^{-n},
 A, D/B, D/C\\
 D, Aq^{1-n}/E, Aq^{1-n}/F \ea \right|\,q,
\;q \right),  \nonumber
\eea
where $ABCq^{1-n} = DEF$. Now apply (4.11) with the invariant parameters
 $A$ and $D$ being the terms depending on $z$. After some routine calculations 
we obtain
\bea
&& \frac{(qb\a, qb\b; q)_{n-1}}
{(q^{1-n}z/\a, q^{1-n}z/\b; q)_n}
\left( \frac{z}{\a\b b q^{n-1}} \right)^{n-1} \frac{z^2}{\a\b } 
 L_n(z) \nonumber   \\
&& = q^{n-1} \sum_{k=0}^{n-1}\frac{(q^{1-n}, az,
 bz, q^{-n+1/2}z/\a\b; q)_k}
 {(q, q^{1-n}z/\a, q^{1-n}z/\b, ab q^{1/2}z; q)_k} q^k \nonumber \\
&& \mbox{} \quad - \frac{(\a\b - q^{1/2}z)} {(\a\b - q^{-n+1/2}z)}
\;\sum_{k=0}^{n-1} \frac{(q^{1-n};q)_k (az, bz, 
q^{-n+1/2}z/\a\b; q)_{k+1}}
 {(q;q)_k(q^{1-n}z/\a, q^{1-n}z/\b, ab q^{1/2}z;q)_{k+1}} q^k. \nonumber
\eea
In the first  sum in the above equation isolate the  $k = 0$ term from the remaining
terms but in the second sum isolate the $k = n-1$ term from the rest
 of the series, then combine the two remaining series. A fantastic cancellation occurs and we obtain
\bea
&&T_{q,z}\left[(\a\b q^{1/2}/z;q)_2 w_c(z; a, q\a, b, q\b|q)
r_{n-1}(z;  a, q\a, b, q\b|q)\right] \\
&& \mbox{} \qquad =\frac{(1-b\a)(1-b\b/q)}
{(1-q)b} w_c(z; a, \a, b, \b|q)
r_n(z; a, \a, b, \b|q), \nonumber
\eea
which is equivalent to (4.8) and the proof is complete.

\bigskip

\setcounter{equation}{0}

\setcounter{section}{5}

\setcounter{thm}{0}

{5. \bf The Biorthogonality Relation}. Let us denote 
\bea
&{}& I_{m,n}(a, \a, b, \b):= \\
&{}& \mbox{} \qquad \mbox{} \quad \frac{1}{2\pi i} \int_{|z|=1} w_{c}(z; a, \a, b, \b|q)r_n(z; a, \a, b, \b|q) \overline{s_m(z; a, \a, b, \b|q)} \frac{dz}{z}, \nonumber
\eea
which is just the integral on the left-hand side of (1.9). Observing that 
\bea
<T_{q,z}f(z), g(z)>_c \; = \; <f(z), D_{q,z} g(z)>_c \nonumber
\eea
we use (4.12) to write 
\bea
&&\frac{(1-\a b)(1-\b b)}{b(1-q)}I_{m,n}(a, \a, b, \b) \\ 
&&\mbox{} \quad = \frac{1}{2\pi i} \int_{|z|=1}(\a\b q^{1/2}/z;q)_2 w_{c}(z; a, q\a, b, q\b|q) r_{n-1}(z; a, q\a, b, q\b|q) \nonumber \\
&& \mbox{} \qquad \times\overline{D_{q, z} s_m(z; a, \a, b, \b|q)} \; \frac{dz}{z}
\nonumber \\
&&=\b q^{1-m}\frac{(1-\a q^{1/2})(1-\b q^{1/2})(1-q^m)(1-ab\a\b q^{m-1})}{(1-q)(1-a\b)(1-b\b)} \nonumber \\
&&\mbox{} \quad \times  \frac{1}{2\pi i} \int_{|z|=1} w_{c}(z; a, q\a, b, q\b|q) r_{n-1}(z; a, q\a, b, q\b|q) \overline{s_{m-1}(z; a, q\a, b, q\b|q)} \frac{dz}{z},
\nonumber
\eea
by (1.7), (1.8) and (4.1). Hence 
\bea
I_{m,n}(a, \a, b, \b) 
&=& - b\b q\frac{(1-q^{1/2}\a)(1- q^{1/2}\b)(1- q^{-m}) (1-a\a b \b q^{m-1})}{(1-\a b)(1- a\b)(1-b\b)^2} \\
&{}& \mbox{} \quad \times I_{m-1, n-1}(a, q\a , b, q\b).
 \nonumber
\eea
Suppose $m \ge n$. Then iterating (5.3) $n-1$ times we obtain 
\bea
I_{m,n}(a, \a, b, \b) &=& (- b\b)^n q^{{{n+1}\choose{2}}}\frac{(q^{1/2}\a, 
q^{1/2}\b,
q^{-m}, a\a b \b q^{m-1}; q)_n}{(\a b,  a\b, b\b, b\b; q)_n}\\
&{}& \mbox{} \qquad \times I_{m-n,0}(a, q^n \a, b, q^n\b). \nonumber
\eea
If $m > n$ then the use of the identity 
\bea
&{}& T_{q, \overline{z}}\left[(abq^{1/2}z;q)_2 w(\overline{z}; \a, qa, \b, qb|q)
\overline{s_{k-1}(z; qa, \a, qb, \b|q)}\right] \\
&{}& \mbox{} \qquad \mbox{} \qquad \mbox{} \qquad = \frac{(1-a\b)(1-b\b)}{ab\b(1-q)}
 w(z; a, \a, b, \b|q)
\overline{s_{k}(z; a, \a, b, \b|q)} \nonumber
\eea
gives
\bea
&& w(z; a, q^n\a, b, q^n\b|q) \overline{s_{m-n}(z; a, q^n\a, b, q^n\b|q)}  \\
&&\mbox{} \quad = \frac{(ab\b q^n(1-q)}{(1-a\b q^n)(1-b\b q^n)} \nonumber \\
&& \mbox{} \quad \times T_{q,\overline{z}}\left[(q^{1/2}abz;q)_2w(\overline{z};
 q^n\a, qa, q^n\b, qb|q)\overline{s_{m-n}(z; qa, q^n\a, qb, q^n\b|q)}\right].
 \nonumber
\eea
But the integral over the expression on the right-hand side of (5.6)  over 
the unit circle vanishes. This proves the orthogonality when $m > n$. If $m=n$ 
then we have 
\bea
I_{0,0}(a, q^n\a, b, q^n\b) &=&  \frac{1}{2\pi i}\int_{|z| = 1}
w(z;  a, q^n\a, b, q^n\b|q\; \frac{dz}{z} \\
&=& \kappa(a, \a, b, \a)\frac{(a\a, b\a, a\b, b\b;q)_n)}{(q^{1/2}\a,
 q^{1/2}\b;q)_n(ab\a\b;q)_{2n}}, \nonumber
\eea
by (1.11). On the other hand, if $m < n$ then, iterating (5.3) $m-1$ 
times we obtain 
\bea
I_{m,n}(a, \a, b, \b) &=& (-b\b)^m
\frac{(q^{1/2}\a,  q^{1/2}\b, q^{-m}, ab\a\b q^{m-1};q)_m}
{(\a b, a\b, b\b, b\b;q)_m} q^{{{m+1}\choose{2}}} \\
&&{} \mbox{} \quad \times I_{0, n-m}(a, q^m\a, b, q^m\b) \nonumber\\
&=& (-b\b)^m
\frac{(q^{1/2}\a,  q^{1/2}\b, q^{-m}, ab\a\b q^{m-1};q)_m}
{(\a b, a\b, b\b, b\b;q)_m} q^{{{m+1}\choose{2}}}\nonumber \\
&&{} \mbox{} \quad \times \frac{1}{2\pi i} \int_{|z| = 1} w(z; a, q^m\a, b, q^m\b|q)
r_{n-m}(z; a, q^m\a, b, q^m\b|q) \; \frac{dz}{z}. \nonumber
\eea
However, by (4.12), 
\bea
&&w(z; a, q^m\a, b, q^m\b|q)r_{n-m}(z; a, q^m\a, b, q^m\b|q) \\
&=& \mbox{} \quad \frac{b(1-q)}{(1-q^m \a b)(1-q^m b \b)} \nonumber \\
&& \times T_{q, z}\left[(\a \b q^{2m+1/2}/z; q)_2
w(z; a, q^{m+1}\a, b, q^{m+1}\b|q)r_{n-m-1}(z; a, q^{m+1}\a, b, q^{m+1}\b|q)\right]
\nonumber
\eea
and the integral of the right-hand side over the unit circle is zero. So, $I_{m,n}(a, \a, b, \b) = 0$ when $m < n$. This together with (5.6) and (5.7), constitutes an alternate proof of the biorthogonality relation (1.9), which is based on the use of the raising and lowering operators $L^{\pm}$.

Department of Mathematics, University of South Florida, Tampa, Florida 33620

Department of Mathematics, Carleton University, Ottawa, Ontario, Canada K1S 5B6 
\end{document}